\documentclass[11pt]{amsart}

\usepackage{latexsym}
\usepackage{amsmath}
\usepackage{amssymb}
\usepackage{epsfig}
\usepackage{amsfonts}

\newcommand\risS[6]{\raisebox{#1pt}[#5pt][#6pt]{\begin{picture}(#4,15)(0,0)
  \put(0,0){\includegraphics[width=#4pt]{#2.eps}} #3
     \end{picture}}}
\newcommand\smf[1]{\risS{-4}{#1}{}{15}{0}{0}}
\newcommand\laf[1]{\begin{picture}(65,40)(0,0)
          \put(0,0){\includegraphics[width=65pt]{#1.eps}}
                   \end{picture}}
\def\kb#1{[ #1 ]}
\def\wt#1{\widetilde{#1}}
\def\a{\alpha}
\def\b{\beta}
\def\bc{\mathrm{bc}}

\def\wh{\widehat}
\def\wo{\overline}
\def\bu{\bullet}

\def\rr{\mathbb R}

\def\Ga{\Gamma}

\def\de{\delta}

\def\al{\alpha}
\def\be{\beta}

\def\ve{\varepsilon}

\def\vp{\varphi}

\def\cS{\mathcal S}

\def\cF{\mathcal F}

\def\ssu{\subset}

\def\wt{\widetilde}
\def\<{\langle}
\def\>{\rangle}

\newtheorem{thm}{Theorem}[section]

\newtheorem{cor}[thm]{Corollary}

\newtheorem{defn}[thm]{Definition}

\newtheorem{Example}[thm]{Example}
\newtheorem{Remark}[thm]{Remark}

\begin{document}

\title[The Kauffman bracket and the Bollob\'as-Riordan polynomial]%
    {The Kauffman bracket and the Bollob\'as-Riordan polynomial of
    ribbon graphs}
\author[Sergei~Chmutov]{Sergei~Chmutov}
\author[Igor~Pak]{Igor~Pak}
\date{21 May, 2004}

\keywords{Knot invariants, Jones polynomial, Kauffman bracket,
Tutte polynomial, \\
\text{\hskip.42cm Bollob\'as-Riordan} polynomial, ribbon graph}

\begin{abstract}
\noindent
For a ribbon graph~$G$ we consider an alternating link~$L_G$
in the 3-manifold $G\times I$ represented as the product of
the oriented surface $G$ and the unit interval~$I$.
We show that the Kauffman bracket~$\kb{L_G}$ is an evaluation
of the recently introduced Bollob\'as-Riordan polynomial~$R_G$.
This results generalizes the celebrated relation between
Kauffman bracket and Tutte polynomial of planar graphs.
\end{abstract}

\maketitle

\section*{Introduction} \label{s:intro}

The study of polynomial invariants was revolutionized
after Jones' celebrated discovery of Jones polynomials
of knots and links in~$\rr^3$.  In~\cite{Ka1}
L.~Kauffman showed how to avoid Jones' somewhat involved
algebraic tools, and presented a simple combinatorial
definition of a more general polynomial, now called
the Kauffman bracket, associated to (signed) planar graphs,
of which Jones polynomial is an evaluation.  In the same
paper Kauffman also showed that this polynomial is
in fact equivalent to the (signed) Tutte polynomial
well studied in combinatorics literature.

The study of the Kauffman bracket, the Jones polynomial,
and other link polynomials was soon extended to links
in~$G\times I$, the products of a two-dimensional surface~$G$
and the interval~$I=[0,1]$; we refer to~\cite{CR,HP,IK,Tu}
and a survey~\cite{Pr} for various approaches to the subject.  
In recent
years the Kauffman bracket of links in~$G\times I$
received additional attention motivated by new results on
crossing numbers~\cite{AFLT} and the study of virtual
knots~\cite{DK,Man}.  All these papers utilize and develop
the first half of ``Kauffman's program'', that is they
extend Kauffman's combinatorial approach to obtain the
generalized Jones and other polynomial invariants of knots
and links in~$G\times I$.

In this paper we are concerned with the second half
of the ``Kauffman's program'', on establishing the connection
between the Kauffman bracket and the Tutte polynomial.  We show
that the Kauffman bracket of links in~$G\times I$ is an
evaluation of the recently introduced Bollob\'as-Riordan
polynomials of ribbon graphs.  The latter generalize the
classical (dichromatic) Tutte polynomial from general
graphs to graphs on surfaces, and our results is an
extension of Kauffman's relation.  Interestingly, the
Bollob\'as-Riordan polynomials were also introduced with
(very different) knot theoretic applications in
mind~\cite{BR2,BR3}.

The paper is structured as follows.  In the first two
sections we recall definitions of Kauffman's bracket
and the Bollob\'as-Riordan polynomial of ribbon graphs.
In section~3 we construct `medial ribbon graphs' and
state the Main Theorem.  As often appears in these cases,
the proofs of results on (generalizations of) the Tutte
polynomial are quite straightforward, so the proof of the
Main Theorem is postponed till section~5.  In section~4
we extend our results to signed ribbon graphs and derive
the Jones polynomial of links in $G \times I$ as an
appropriate evaluation.  We conclude with final remarks
and overview of the literature.

\bigskip

\section{The Kauffman bracket in $G\times I$.}
Let $G$ be an oriented surface (possibly with the boundary),
let~$I=[0,1]$ be the unit interval,
and let~$L$ be an unoriented link in the 3-manifold~$G\times I$.
To represent~$L$ by its diagram on the surface~$G$ we will always
assume that~$L$ is in general position with respect to the
projection~$\pi\colon G\times I \to G$, i.e. the image~$\pi(L)$
is an immersed curve in~$G$ with finitely many double points
as its only singularities. This can be achieved by a small
deformation of~$L$ which preserves the topological type of the link~$L$.
Then a {\it diagram} of~$L$, denotes~$\wt L$ is the curve
$\pi(L)\subset G$ with an
extra information of {\it overcrossing} and {\it undercrossing}
at every double point of~$\pi(L)$. Naturally, the link~$L$ can have
many different diagrams~$\wt L$.

Consider two ways of resolving a crossing in~$L$.
The {\it $A$-splitting}, $\smf{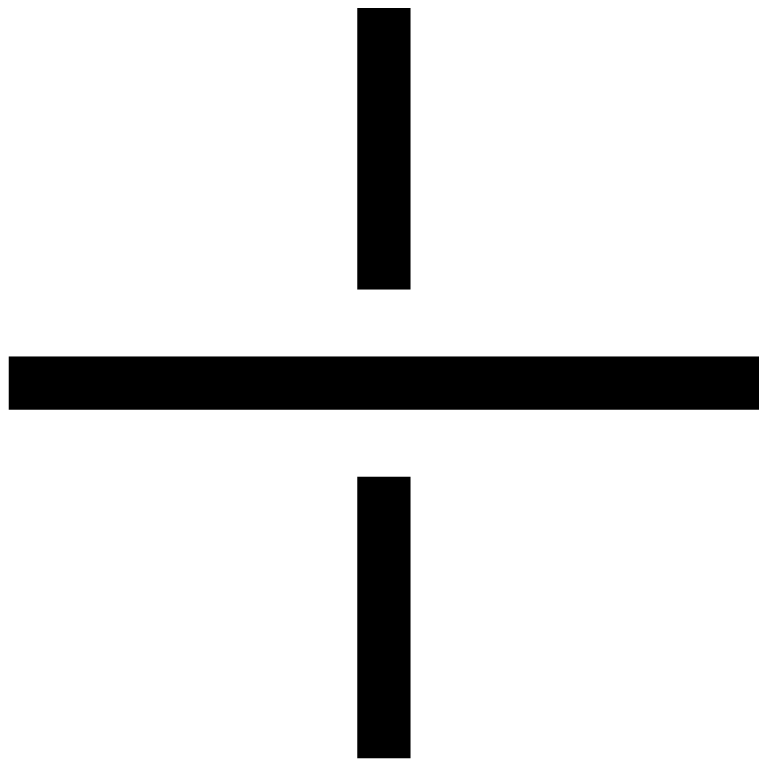}\ \leadsto\ \smf{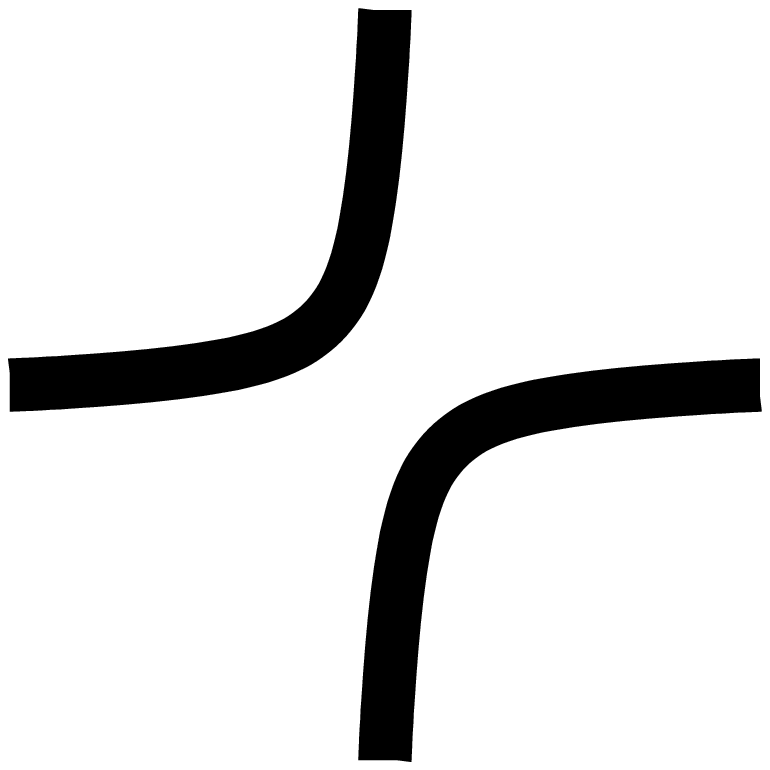}$,
is obtained by uniting two regions swept out by the overcrossing arc under
the   counterclockwise rotation until the undercrossing arc. We are
assuming here that the orientation of~$G$ is given by the counterclockwise
rotation, so the overcrossing arc was rotated according to the orientation.
Similarly, the {\it $B$-splitting},
$\smf{cr}\ \leadsto\ \smf{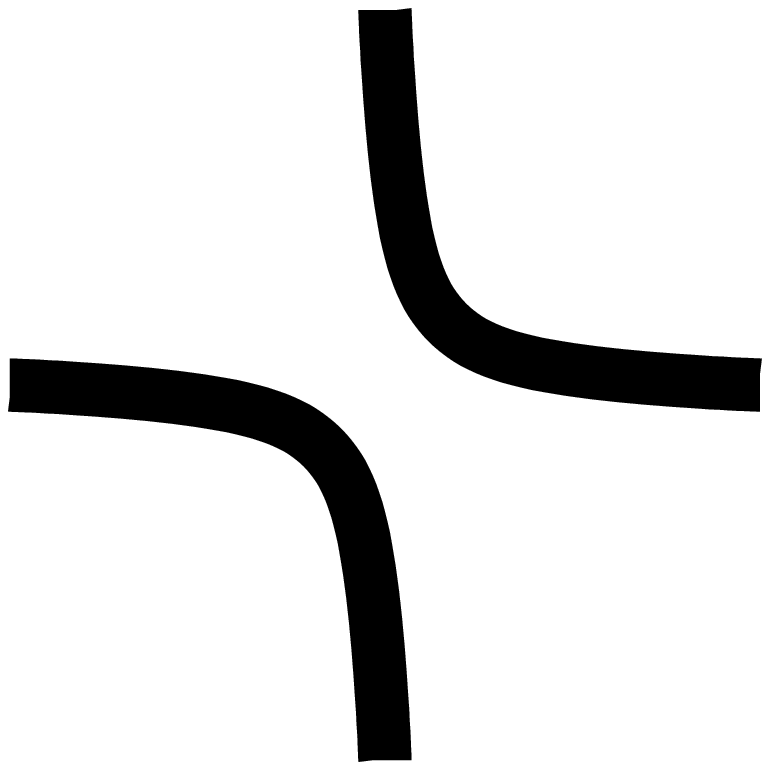}$, is
obtained by uniting two other regions. A {\it state}~$S$ of
the link diagram~$\wt L$
is a way of splitting at each crossing of the diagram, and
denote by~$\cS(\wt L)$ the set of such states.
Clearly, a diagram~$\wt L$ with~$n$ crossings has~$|\cS(\wt L)| = 2^n$
different states.

Denote by~$\a(S)$ and~$\b(S)$ the number of $A$-splittings and $B$-splittings
in a state~$S$, respectively.  Also, denote by~$\de(S)$ the number of
connected components of the curve obtained from the link
diagram~$\wt L$ by all
splittings according to state~$S \in \cS(\wt L)$.

\begin{defn}\label{def:kb}
The \emph{Kauffman bracket} of a diagram $\wt L$ of a
link $L\subset G\times I$ is a polynomial in three variables
$A$, $B$, $d$ defined by the formula:
\begin{equation}
\kb{\wt L} (A,B,d)\ :=\ \sum_{S \in \cS(\wt L)} \,
A^{\a(S)} \, B^{\b(S)} \, d^{\de(S)-1}\,.
\end{equation}
\end{defn}

This definition of the (generalized) Kauffman's bracket
follows~\cite{Kam}.
Note that $\kb{\wt L}$ is \emph{not} a topological
invariant of the link~$L$ and in fact depends on the link
diagram~$\wt L$.

\begin{Example}\label{ex1} {\rm
Consider the surface~$G$ and the link
$L\subset G\times I$ as shown on the left in the table below.
The corresponding diagram~$\wt L$ has two crossings, so there are
four states for it,
$|\cS(\wt L)| = 4$.
The curves obtained by the splittings and the corresponding
parameters $\a(S)$, $\b(S)$, and $\de(S)$ are shown in the
remaining columns of the table.
$$
\begin{array}{c||c|c|c|c}
\raisebox{-8pt}{\begin{picture}(80,75)(0,-3)
\put(0,0){\includegraphics[width=80pt]{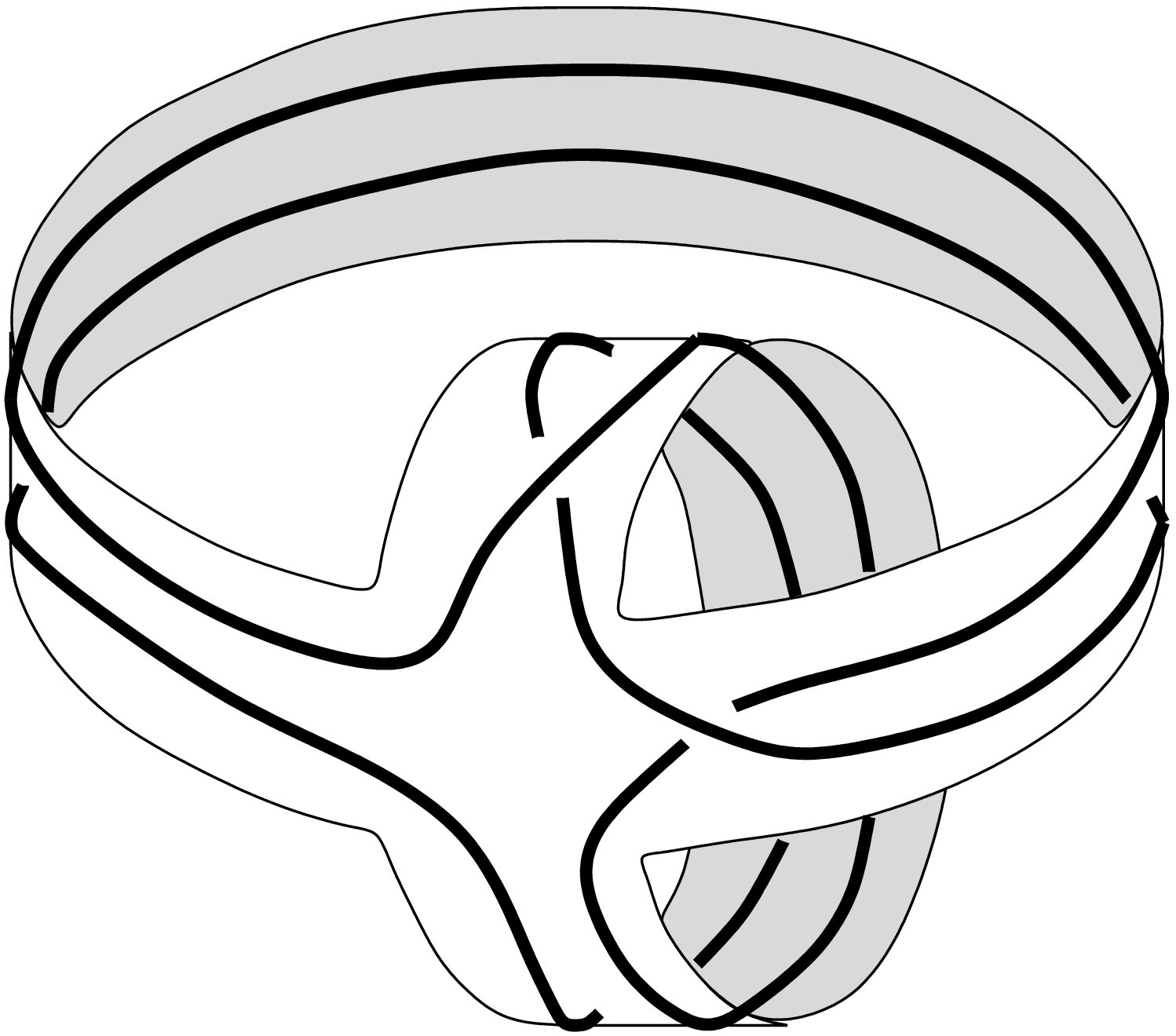}}\end{picture}\vspace*{10pt}}
 & \laf{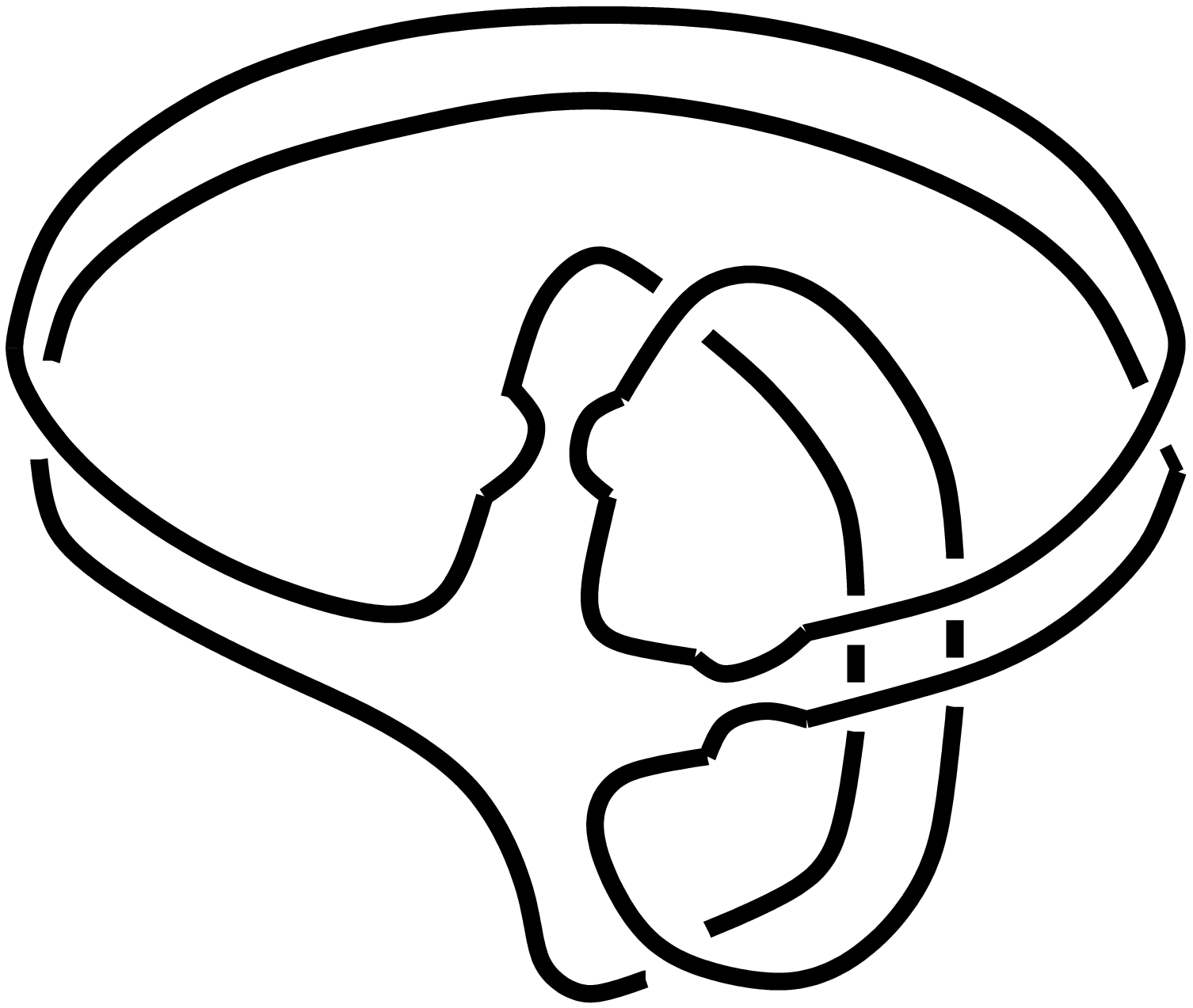} & \laf{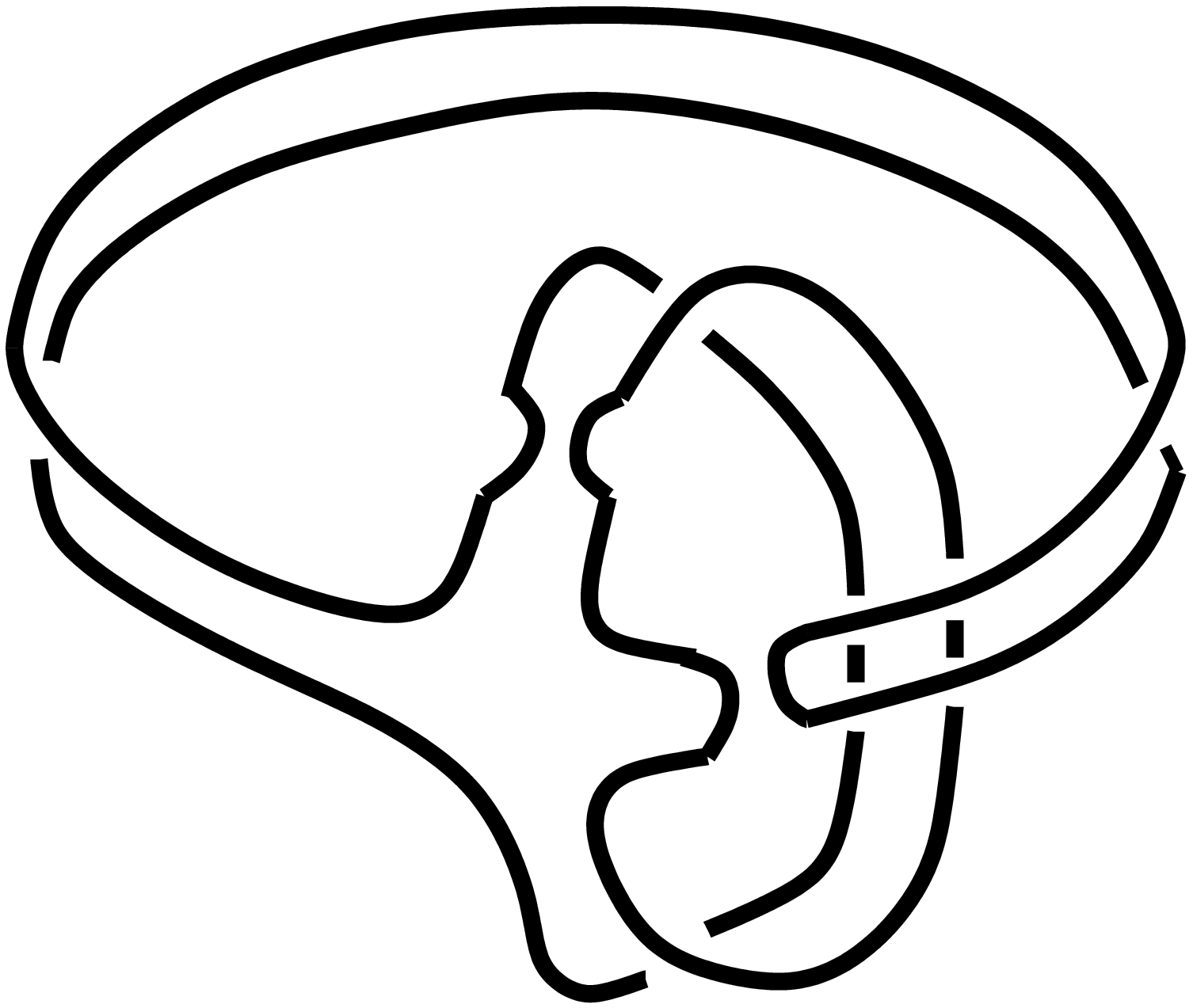} & \laf{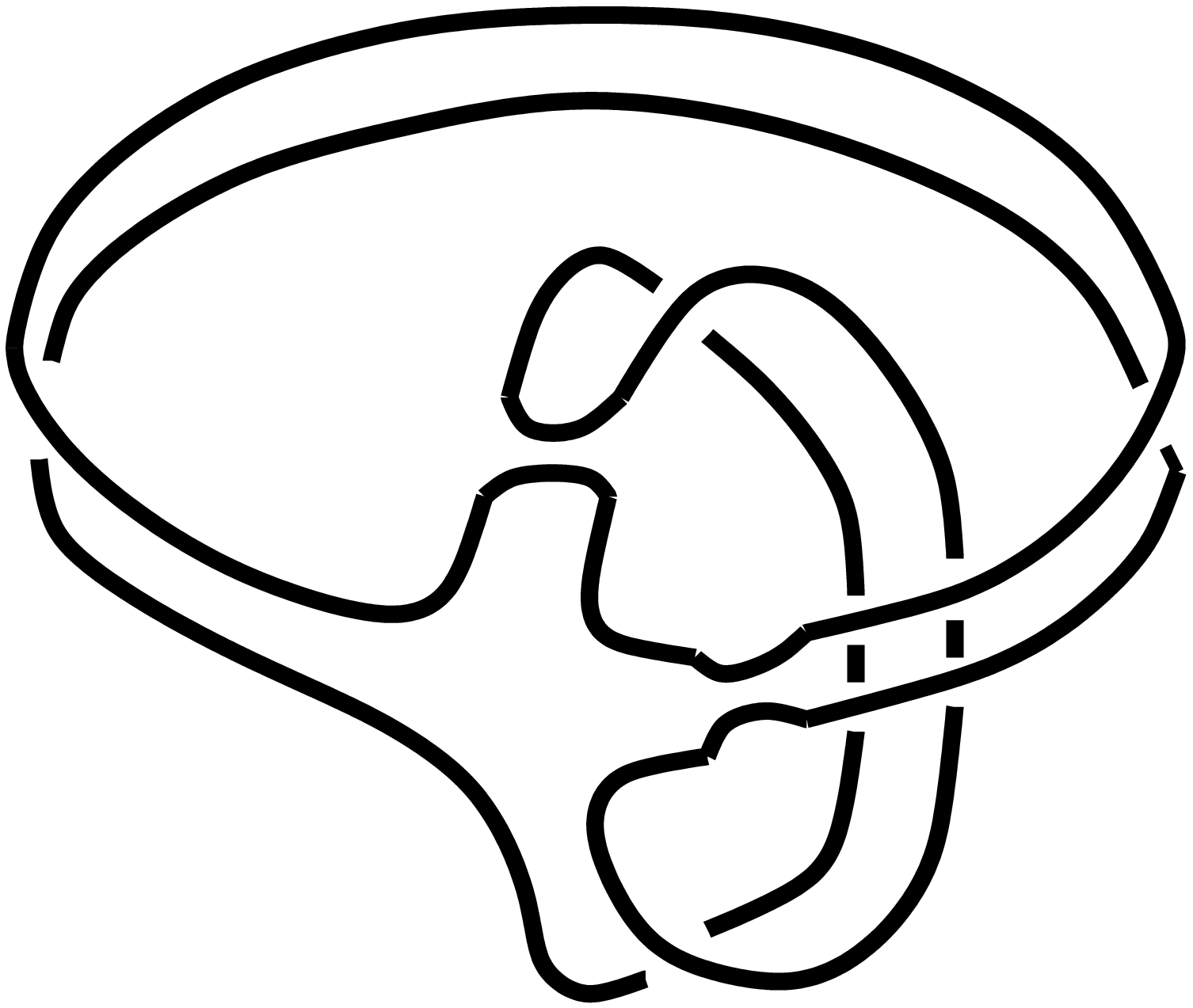} & \laf{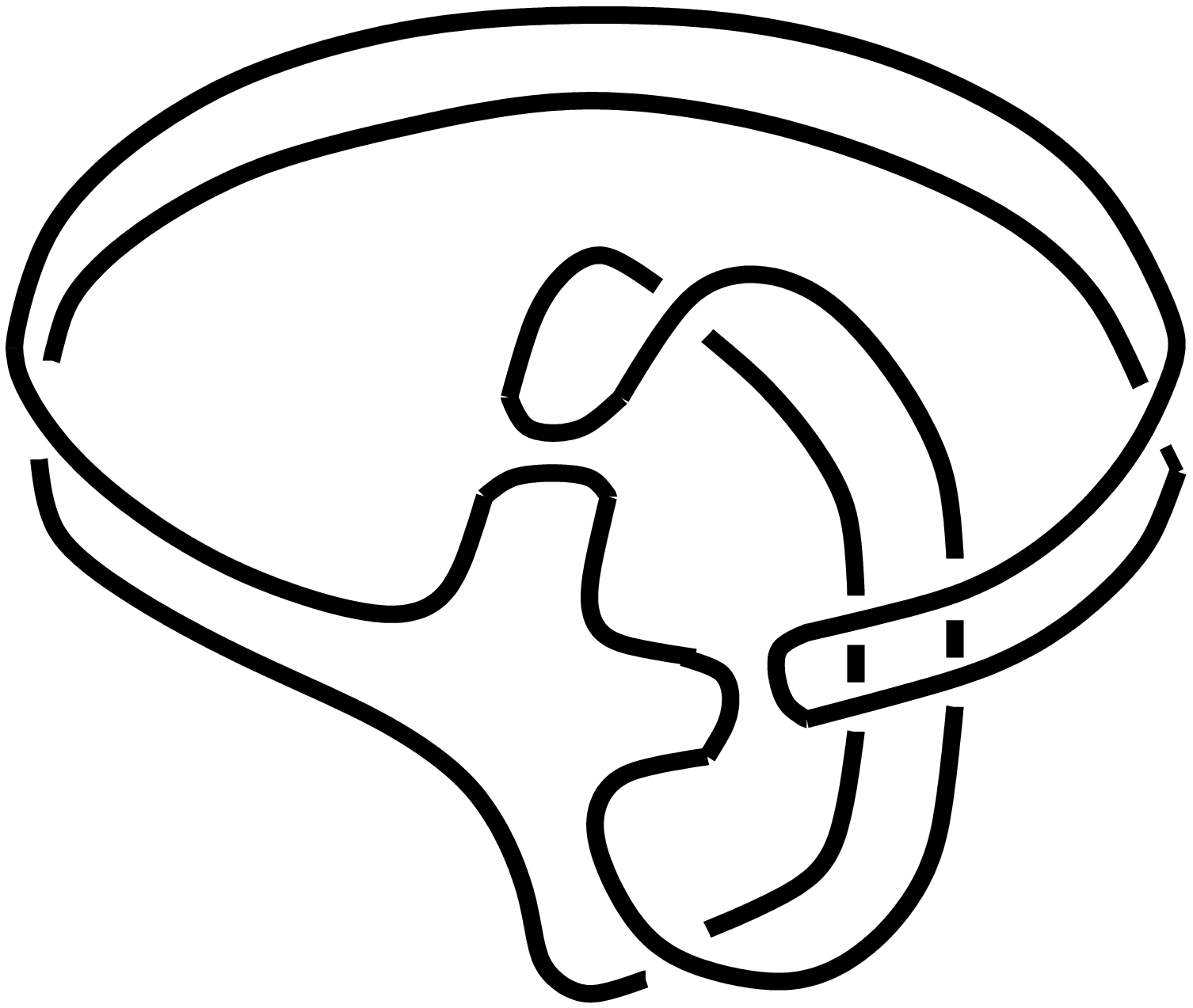}\\ \hline
 (\a(S),\b(S),\de(S)) & (2,0,1) & (1,1,2) & (1,1,2) & (0,2,1)\makebox(0,15){}
\end{array}
$$
In this case the Kauffman bracket of~$\wt L$ is given by \,
$\kb{\wt L} = A^2+ 2 A B d  + B^2.$
}
\end{Example}

\bigskip

\section{The Bollob\'as-Riordan polynomial.}

Let~$\Ga = (V,E)$ be a undirected graph with the set of vertices~$V$
and the set of edges~$E$ (loops and multiple edges are allowed).
Suppose in each vertex~$v \in V$ there is a fixed cyclic order on
edges adjacent to~$v$ (loops are counted twice).  We call this
combinatorial structure the {\it ribbon graph}, and denote it
by~$G$.  One can represent~$G$ by making vertices into `discs'
and connecting them by `ribbons' as prescribed by the cyclic
orders (see Example~\ref{ex2} below).  This defines a 2-dimensional
surface with the boundary, which by a
slight abuse of notation we also denote by~$G$.

Formally, $G$ is the surface with the boundary represented as
the union of two sets of closed topological disks, corresponding to
vertices~$v \in V$ and edges~$e \in E$, which satisfies the following
conditions:

\smallskip

$\bu$ \ these discs and ribbons intersect by disjoint line segments,

$\bu$ \ each such line segment lies on the boundary of precisely
one vertex and precisely\\
\hspace*{.75cm}  one edge,

$\bu$ \  every edge contains exactly two such line segments.

\smallskip

\noindent
It will be clear from the context whether by~$G$ we mean the
ribbon graph or its underlying surface.  In this paper we
restrict ourselves to oriented surfaces~$G$. 
We refer to~\cite{GT} for other definitions and references. 
We use the following standard notation. 

For a ribbon graph~$G$ let
$v(G) = |V|$ denotes the number of vertices,
$e(G) = |E|$ denotes the number of edges, and
let~$k(G)$ be the number of connected components of~$G$.
Also, let~$r(G)=v(G)-k(G)$ be the {\it rank} of~$G$, and
let~$n(G)=e(G)-r(G)$ be the {\it nullity} of~$G$.
Finally, let~$\bc(G)$ be the number of connected
components of the boundary of the surface~$G$.

A spanning subgraph of a ribbon graph~$G$ is defined as a
subgraph which contains all the vertices, and a subset of edges.
Let $\cF(G)$ denote the set of spanning subgraphs of~$G$.
Clearly, the number of spanning subgraphs $F\subseteq G$
is equal to~$|\cF(G)| = 2^{e(G)}$.

\begin{defn} \label{def:br}
The \emph{Bollob\'as-Riordan polynomial} $R_G(x,y,z)$ of a ribbon
graph $G$ is defined by the formula
\begin{equation}\label{def_br}
R_G(x,y,z)\ :=\ \sum_{F \in \cF(G)}
   x^{r(G)-r(F)} y^{n(F)} z^{k(F)-\bc(F)+n(F)}\,,
\end{equation}
where the sum runs over all spanning subgraphs $F$ of $G$.
\end{defn}

This version of the polynomial is obtained from the original
one~\cite{BR2,BR3} by a simple substitution.
Note that for all planar ribbon graphs~$F\subseteq G$
(i.e. when the surface $F$ has genus zero)
the Euler formula gives $k(F)-\bc(F)+n(F)=0$. So, for a planar
ribbon graph~$G$ the Bollob\'as-Riordan polynomial $R_G$ does not
contain powers of~$z$. In fact, in this case it is essentially
equal to the classical Tutte polynomial $T_{\Ga}(x,y)$ of the
(abstract) graph~$\Ga$:
$$R_G(x-1,y-1,z) = T_{\Ga}(x,y)\,.$$
Similarly, a specialization $z=1$ of the Bollob\'as-Riordan polynomial
of arbitrary ribbon graph~$G$, gives the Tutte polynomial once
again:
$$R_G(x-1,y-1,1) = T_{\Ga}(x,y)\,.$$
We refer to~\cite{BR2,BR3} for proofs of these formulas and
to~\cite{B,W} for general background on the Tutte polynomial.

\begin{Example}\label{ex2} {\rm
Consider the following ribbon graph~$G$ with the surface as in the
Example~\ref{ex1}, which corresponds to the abstract graph~$\Ga$ as
below.

$$G\ =\quad \raisebox{-40pt}{\begin{picture}(100,85)(0,0)
  \put(0,0){\includegraphics[width=100pt]{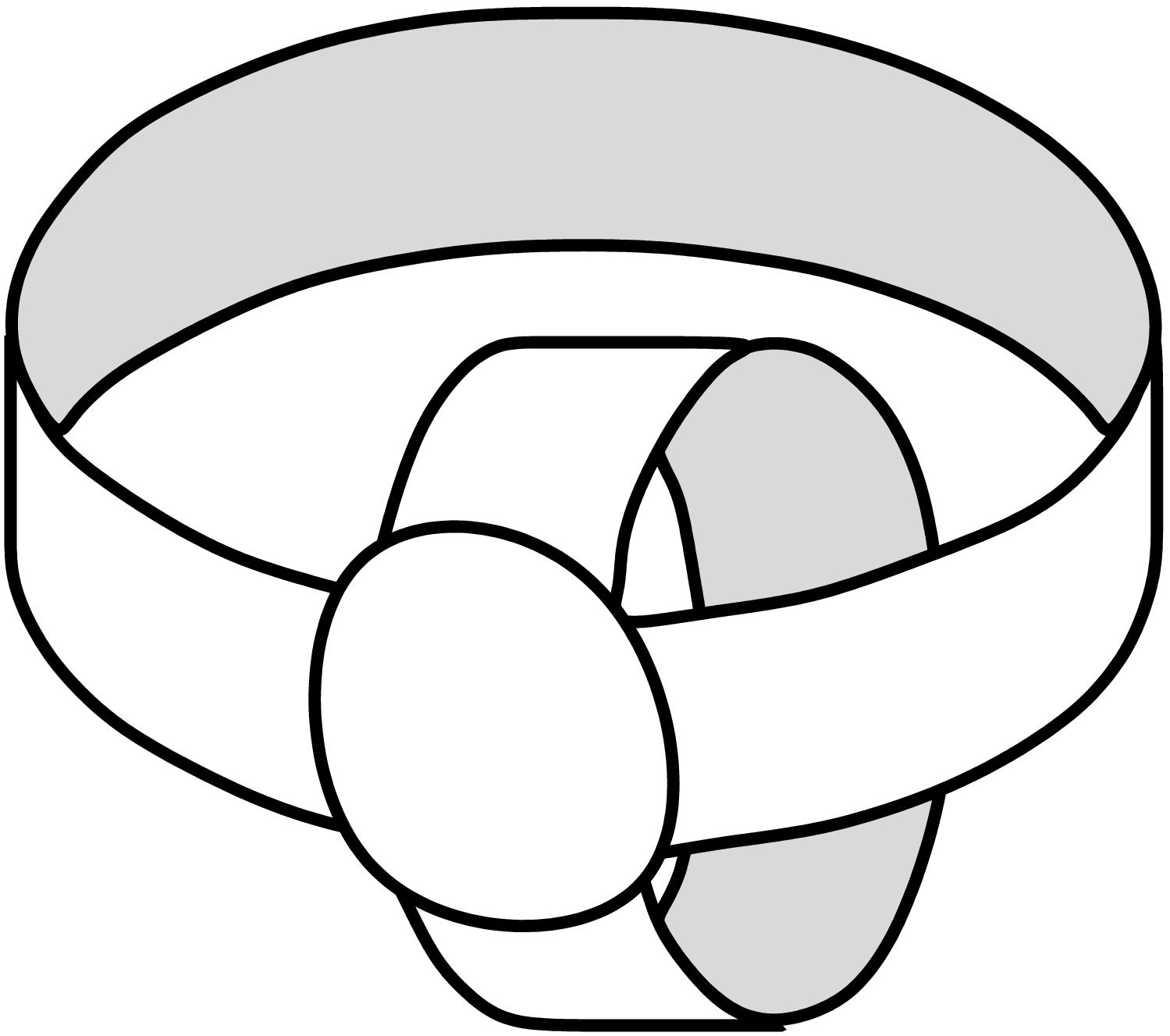}}
     \end{picture} } \qquad\qquad
\begin{array}{rcl}
v(G)&=& 1\\
e(G)&=& 2\\
k(G)&=& 1\\
r(G)&=& 0\\
n(G)&=& 2\\
\bc(G)&=& 1\\
\end{array}\qquad\qquad
\Ga\ =\ \raisebox{-15pt}{\begin{picture}(50,50)(0,0)
  \put(0,0){\includegraphics[width=50pt]{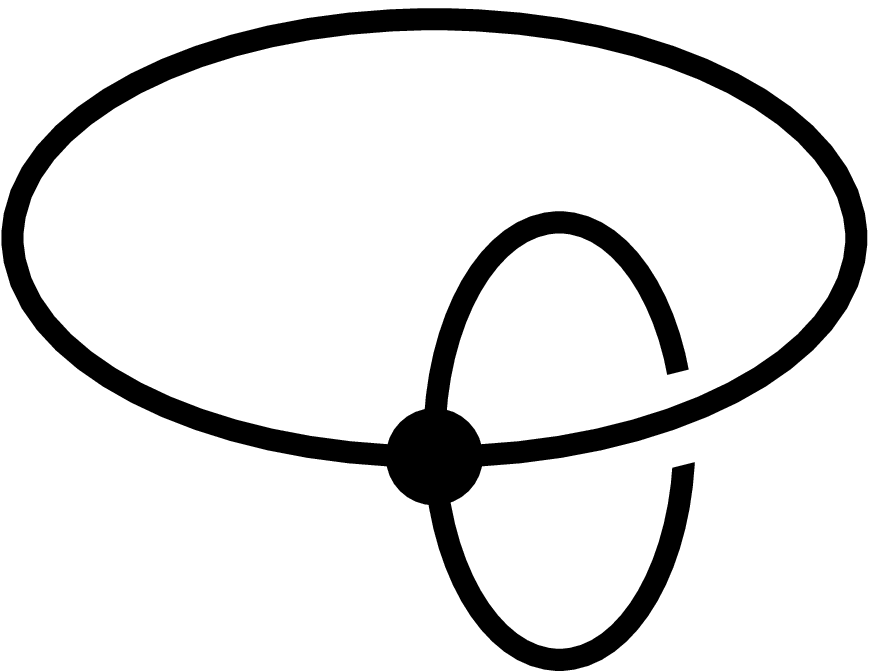}}
     \end{picture} }
$$

By definition, one can compute the corresponding
Bollob\'as-Riordan polynomial:
$$R_G(x,y,z) \, = \, 1+2y+y^2z^2\,.$$
}
\end{Example}

\bigskip

\section{Medial alternating links and the Main Theorem.}

The main result of this paper is a connection between the
Bollob\'as-Riordan polynomial $R_G(x,y,z)$ and the Kauffman
bracket $\kb{\wt L_G}(A,B,d)$ of a medial alternating link
$L_G \in G \times I$ defined below.  This connection
naturally generalizes Kauffman's result for planar
graphs.

\smallskip

Let $G$ be a ribbon graph and $\Ga$ be the corresponding abstract
graph embedded into the surface~$G$.  We construct the (ribbon)
{\it medial}
{\it graph}~$H_{G}$ by embedding into the surface~$G$ as follows.
The vertices of $H_{G}$ lie in the middle of the edges
of~$\Ga$; every vertex has valence~4. The edges of $H_G$
go along the edge-ribbon of~$G$ until the intersection with the
vertex-disks where they turn to the next edge-ribbon
(see the figure below).

Now, consider the \emph{chess-board coloring} of the regions of
$H_G\colon$ color a region \emph{black} if it contains a vertex
of~$\Ga$; otherwise color it \emph{white}.  We define the
link~$L_G\subset G\times I$ as follows. Make a link diagram~$\wt L_G$
out of $H_G$ by making a crossing at each vertex of $H_{G}$
in such a way that the overcrossing branch sweeps out the black regions
in our coloring when rotated according to the orientation of the
surface $G$ until the undercrossing branch. See an example below,
where the orientation of $G$ is given by the
counterclockwise rotation:
$$\ \raisebox{-59pt}{\begin{picture}(120,125)(0,0)
  \put(0,0){\includegraphics[width=120pt]{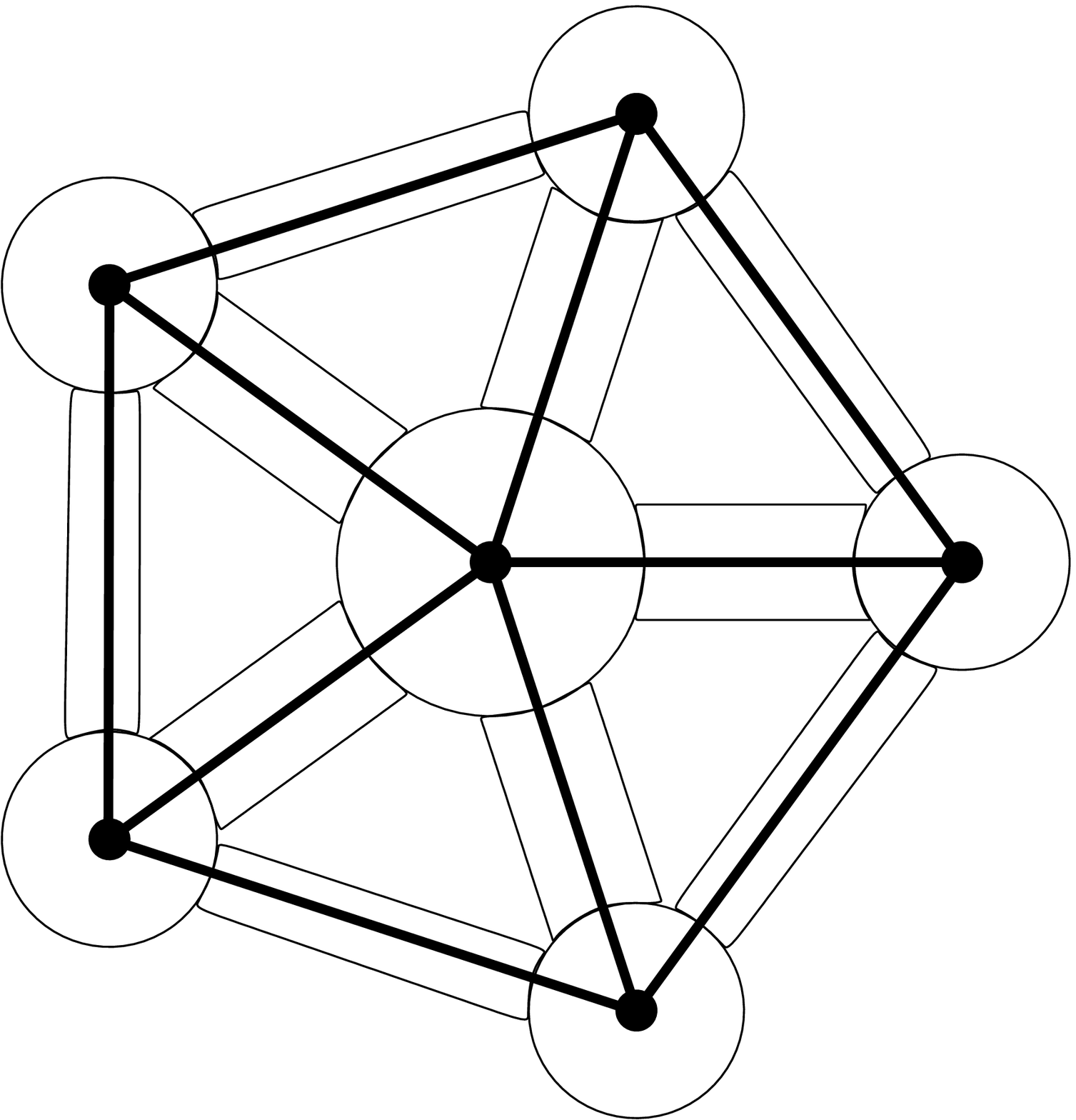}}
  \put(-1,47){\mbox{$\scriptstyle G$}}
  \put(49,50){\mbox{$\scriptstyle \Ga$}}
     \end{picture} }  \leadsto
\raisebox{-59pt}{\begin{picture}(120,110)(0,0)
  \put(0,0){\includegraphics[width=120pt]{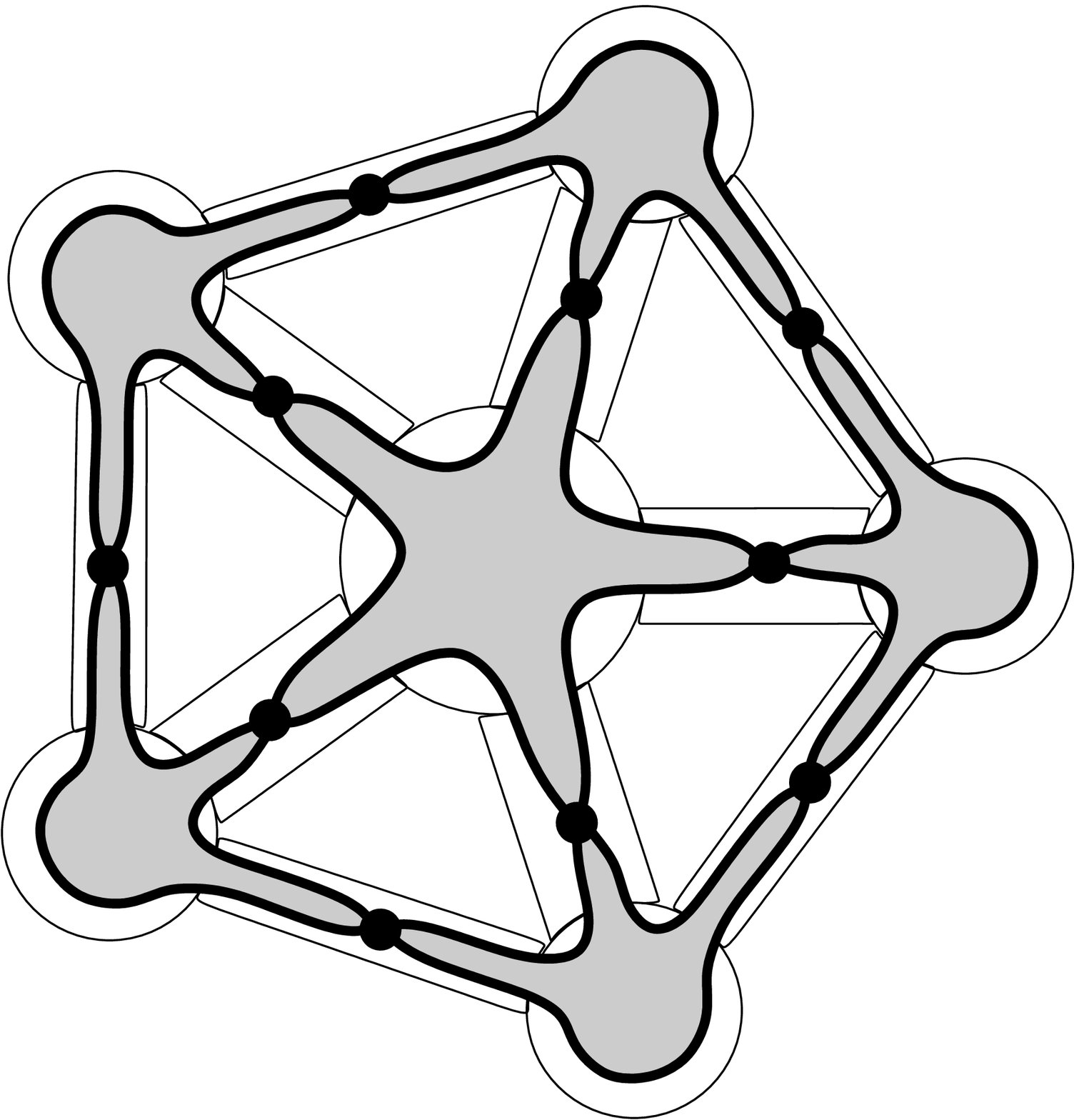}}
  \put(-1,48){\mbox{$\scriptstyle G$}}
  \put(47,60){\mbox{$\scriptstyle H_G$}}
     \end{picture} }  \leadsto
\raisebox{-59pt}{\begin{picture}(120,110)(0,0)
  \put(0,0){\includegraphics[width=120pt]{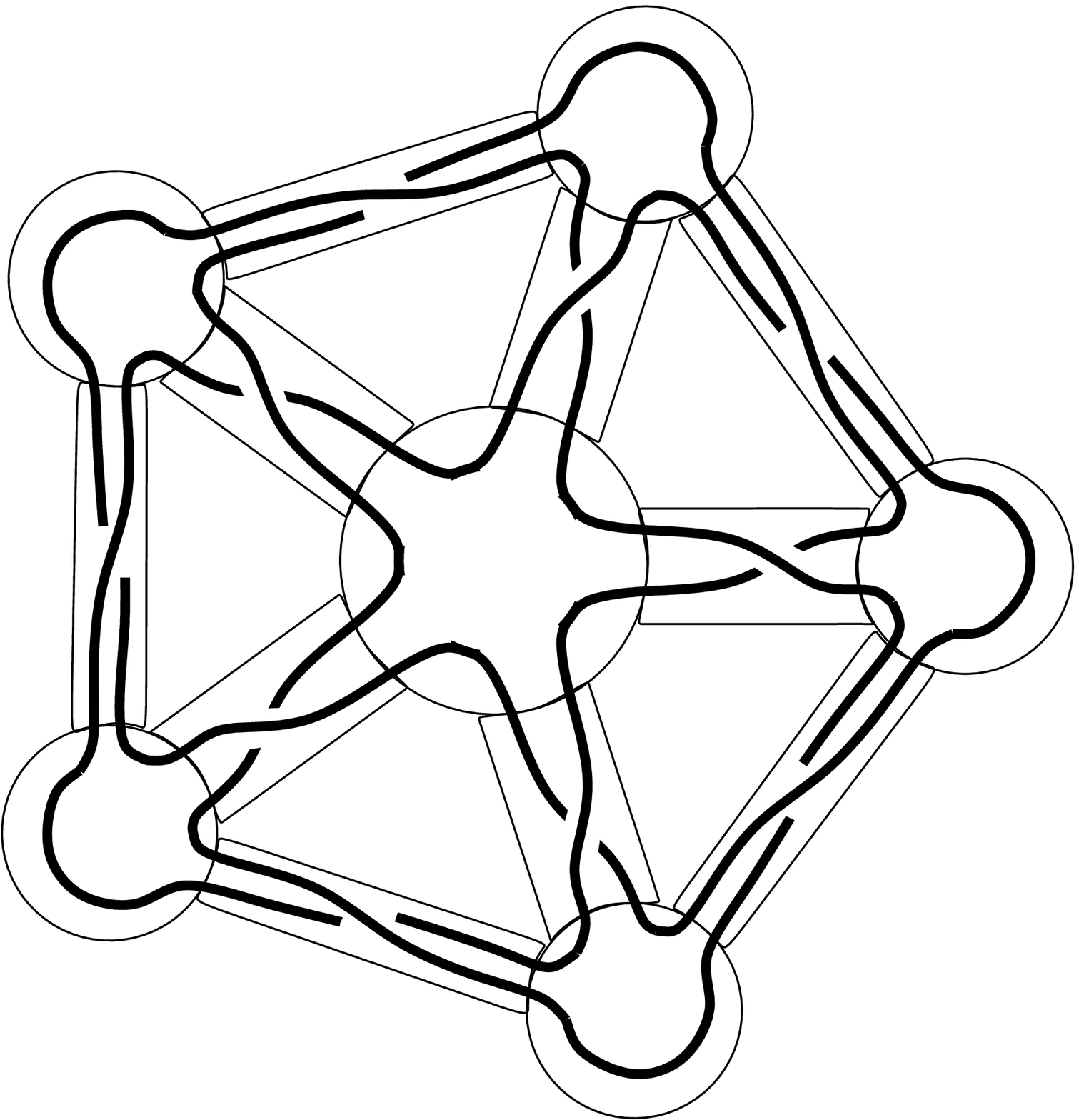}}
  \put(-1,48){\mbox{$\scriptstyle G$}}
  \put(47,60){\mbox{$\scriptstyle \wt L_G$}}
     \end{picture} }
$$

This construction gives the \emph{natural diagram}~$\wt L_G$ of the
{\it medial alternating link $L_G$} for the ribbon graph $G$.
Note that the image of $L_G$ under the projection $G\times I \to G$
is the medial graph~$H_G$. Also, the medial alternating link $L_G$
for the ribbon graph in Example~\ref{ex2}
is precisely the link~$L$ in Example~\ref{ex1}.

\medskip

\noindent
{\bf Main Theorem 3.1.}  {\it For every ribbon graph~$G$ and
the corresponding medial alternating link $L_G \ssu G \times I$
with the natural diagram~$\wt L_G$ we have:
$$\kb{\wt L_G} (A,B,d) \, = \, A^{r(G)} B^{n(G)} d^{k(G)-1} \,
R_G\left(\frac{Bd}{A}, \frac{Ad}{B}, \frac{1}{d}\right).$$
}

\medskip

We should warn the reader that although both Kauffman bracket
and Bollob\'as-Riordan polynomial are polynomials in three
variables, the former has only two free variables since
$\kb{\wt L}$ is always homogeneous in~$A$ and~$B$.
This follows from $\al(S) + \be(S) = e(S)$ for
all~$S \in \cS(\wt L)$.  Another way to see this is to note that
the values of the Bollob\'as-Riordan polynomial $R_G(x,y,z)$
evaluated in the theorem satisfy  $xyz^2 = 1$ equation.
So the situation here is different from the planar case where the
Kauffman bracket $\kb{\wt L_G}(A,B,d)$ and the Tutte polynonomial
$T_{\Ga}(x,y)$ detrmine each other.

\bigskip

\section{Extensions and applications.}

The \emph{signed ribbon graph}~$\wh G$ is a ribbon
graph~$G$ given by~$(V,E)$ with a \emph{sign function}
$\ve: E \to \{\pm 1\}$.
For a spanning subgraph $F \ssu G$ denote by~$e_{-}(F)$
the number of edges~$e \in E$ with~$\ve(e)=-1$.
Denote by~$\wo F = G - F$ the complementary spanning
subgraph of~$G$ with only those edges of~$G$ that
do not belong to~$F$. Finally,
let~$s(F)=\frac12 \bigl(e_{-}(F)-e_{-}(\wo F)\bigr)$.

 We define the
\emph{signed Bollob\'as-Riordan polynomial} $R_{\wh G}(x,y,z)$
as follows:
\begin{equation}\label{def_sbr}
R_{\wh G}(x,y,z)\ :=\ \sum_{F \in \cF(G)}
   x^{r(G)-r(F)+s(F)}
   y^{n(F)-s(F)}
   z^{k(F)-\bc(F)+n(F)}\, .
\end{equation}

\medskip

Similarly, define the \emph{signed medial link} $L_{\wh G}$ by a diagram
obtained from $\wt L_G$ by switching the overcrossings to undercrossings
for negative edges of $\wh G$.

\medskip

\begin{thm} \label{t:signed}
For every signed ribbon graph~$\wh G$ and
the corresponding signed medial link $L_{\wh G} \ssu G \times I$
with the diagram~$\wt L = \wt L_{\wh G}$ we have:
$$\kb{\wt L} (A,B,d)  \ = \ A^{r(G)} B^{n(G)} d^{k(G)-1} \,
R_{\wh G}\left(\frac{Bd}{A}, \frac{Ad}{B}, \frac{1}{d}\right).$$
\end{thm}

\smallskip

The proof is follows verbatim the proof of the Main Theorem
(see the next section).  We leave the details to the reader.

\smallskip

Now recall that the Jones polynomial is defined for oriented
links.  It can be obtained
evaluating the Kauffman bracket~\cite{Ka1}:
$$J_L(t) \, := \,
  (-1)^{w(\wt L)} \, t^{3w(\wt L)/4} \,
\kb{\wt L}(t^{-1/4}, t^{1/4}, -t^{1/2}-t^{-1/2}) \, ,
$$
where $w(\wt L)$ is the \emph{writhe} of a diagram $\wt L$ determined
by the orientation of $\wt L$ as the sum of the following signs of the
crossings of~$\wt L$\,:
$$\begin{picture}(50,30)(0,0)
  \put(0,0){\includegraphics[width=50pt]{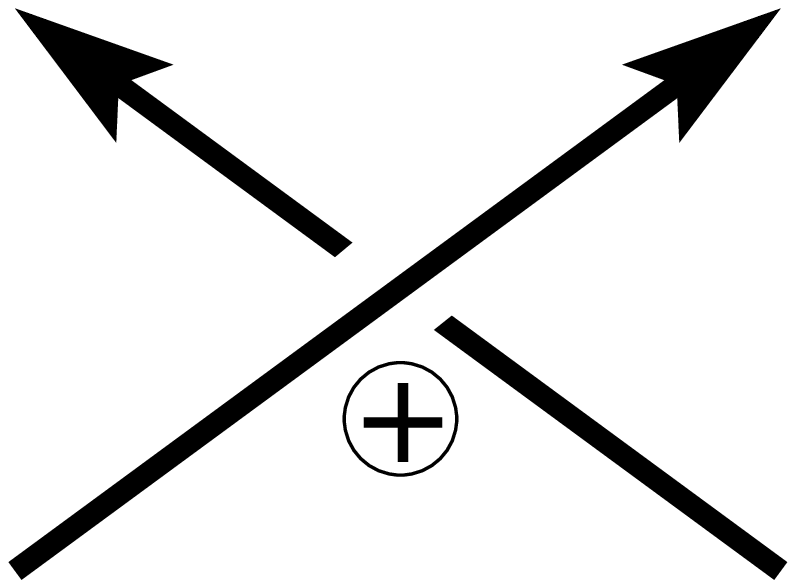}}
     \end{picture} \hspace{3cm}
\begin{picture}(50,40)(0,0)
  \put(0,0){\includegraphics[width=50pt]{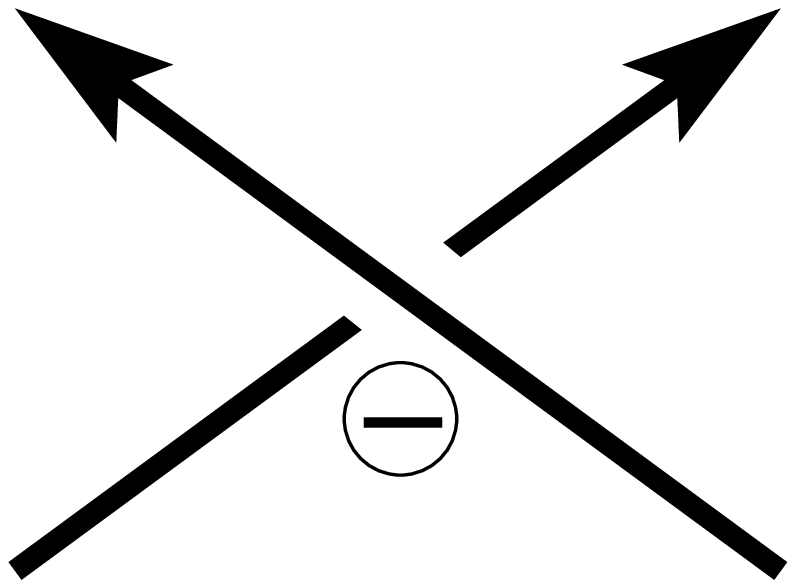}}
     \end{picture}
$$
The following result is an immediate consequence of Theorem~\ref{t:signed}.

\medskip

\begin{cor} \label{c:jones}
Let $\wh G$ be a signed ribbon graph and let
$L=L_{\wh G}$ be the corresponding signed medial link endowed
with an orientation. Then
$$J_L(t) \, = \, (-1)^{w(\wt L)} \,
t^{\frac{3w(\wt L)-r(\wh G)+n(\wh G)}{4}} \,
      \bigl(-t^{1/2}\!-t^{-1/2}\bigr)^{k(\wh G)-1} \,
R_{\wh G}\Bigl(-t-1, -t^{-1}\!-1, \frac{1}{-t^{1/2}\!-t^{-1/2}} \Bigr).
$$
In particular, if $\wh G$ is a planar ribbon graph with only positive edges
and with underlying graph $\Ga$, we have the following well-known relation:
$$J_L(t) \, = \, (-1)^{w(\wt L)} t^{\frac{3w(\wt L)-r(\wh G)+n(\wh G)}{4}}
    \,  \bigl(-t^{1/2}-t^{-1/2}\bigr)^{k(\wh G)-1}  \, T_{\Ga}(-t, -t^{-1}).
$$
\end{cor}

\bigskip

\begin{Remark} \, {\rm
In contrast with the (usual) links in~$\rr^3$,
not every link in the product of a surface~$G$ and
the interval~$I$ can be represented as a signed medial link of
a subsurface of~$G$. The simplest example is the link in the product
of the torus and~$I$, as shown in figure below.
Thus Theorem~\ref{t:signed} and Corollary~\ref{c:jones}
cannot be used to compute the Kauffman bracket and the
Jones polynomial of an arbitrary link~$L \in G \times I$.
\vskip.42cm
\hskip5.cm
\begin{picture}(120,25)(0,57)
\put(0,0){\includegraphics[width=120pt]{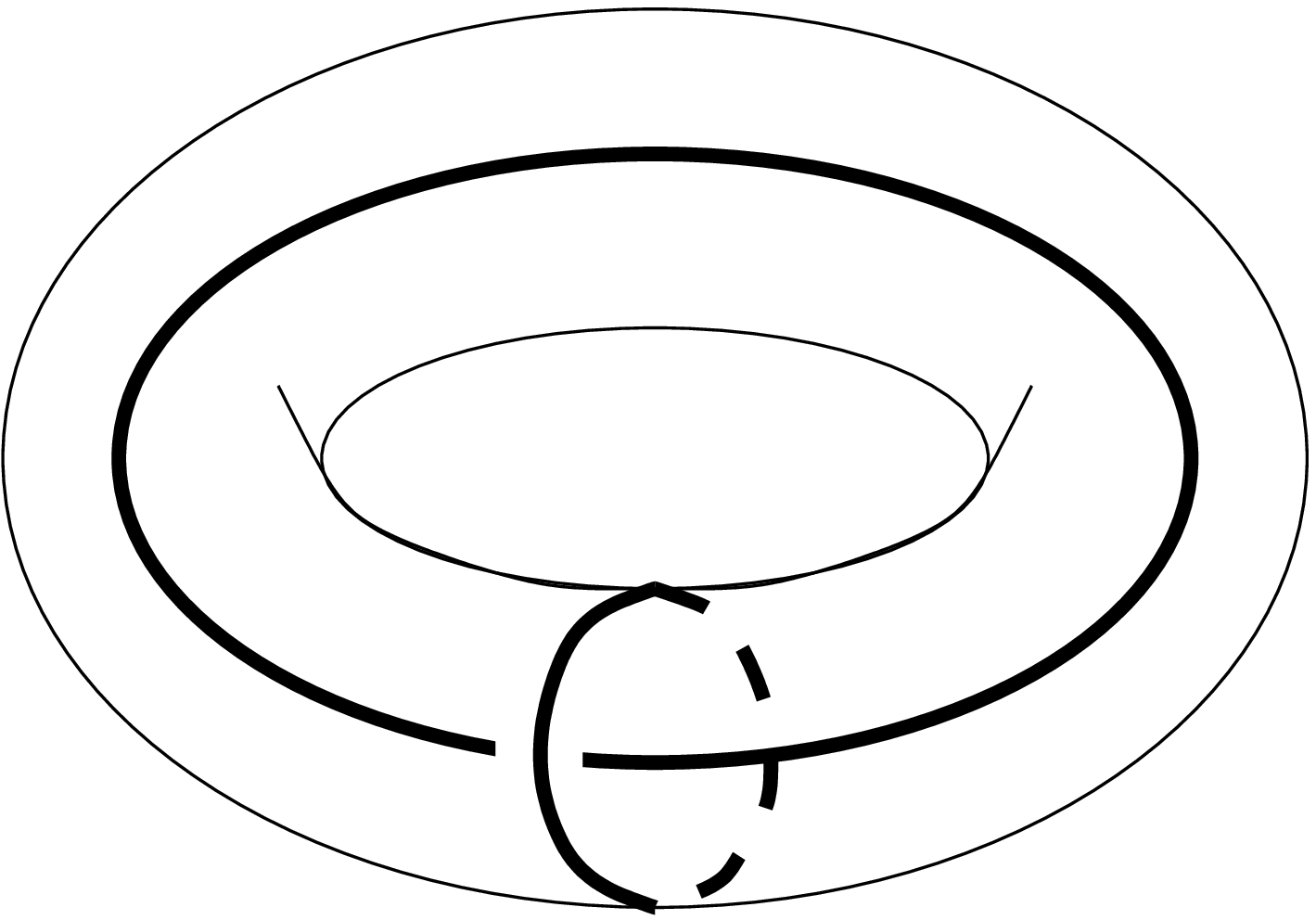}}\end{picture}
}
\end{Remark}

\vskip2.5cm

\section{Proof of the Main Theorem.}
The notation used in the definitions of the Kauffman bracket
(Definition~\ref{def:kb}) and the Bollob\'as-Riordan polynomial
(Definition~\ref{def:br}) suggests a hint on how to prove the
Main Theorem.  There is a
natural one-to-one correspondence $\vp: \cS(\wt L_G) \to \cF(G)$
between the states $S \in \cS(\wt L_G)$ and spanning
subgraphs~$F\subseteq G$.
Indeed, note that the crossings of the diagram~$\wt L_G$
correspond to the edges of~$G$.  Now,
let an $A$-splitting of a crossing in~$S$ mean that we keep the
corresponding edge in the spanning subgraph~$F= \vp(S)$.
Similarly let a $B$-splitting in~$S$ mean that we remove the
edge from the subgraph~$F = \vp(S)$.

By definition, we have~$\de(S) = \bc(F)$, for all~$F = \vp(S)$.
Furthermore, we easily obtain the following
relation between the parameters:
$$e(F)=\a(S),\qquad e(G)-e(F)=\b(S)\ ,$$
for all $S \in \cS(\wt L_G)$, and $F = \vp(S)$.
Now, for a spanning subgraph $F \in \cF(G)$, consider a term
$x^{r(G)-r(F)} y^{n(F)} z^{k(F)-\bc(F)+n(F)}$ of $R_G(x,y,z)$.
After the substitution
$$x\, = \, \frac{Bd}{A},\qquad y \, = \, \frac{Ad}{B},\qquad z\, =\,
\frac{1}{d}$$
and multiplication of this term by
$A^{r(G)} B^{n(G)} d^{k(G)-1}$ as in the Main Theorem,
we get
$$\begin{array}{ccl}
\multicolumn{3}{l}{ A^{r(G)} B^{n(G)} d^{k(G)-1}
                    (A^{-1}Bd)^{r(G)-r(F)} (AB^{-1}d)^{n(F)}
                    d^{-k(F)+\bc(F)-n(F)}  }     \vspace{12pt}\\
&=\!& A^{r(G)-r(G)+r(F)+n(F)} B^{n(G)+r(G)-r(F)-n(F)}
          d^{k(G)-1+r(G)-r(F)+n(F)-k(F)+\bc(F)-n(F)}  \vspace{12pt}\\
&=\!& A^{r(F)+n(F)} B^{n(G)+r(G)-r(F)-n(F)}
          d^{k(G)-1+r(G)-r(F)-k(F)+\bc(F)}\ .
\end{array}$$

It is easy to see that $r(F)+n(F)=e(F)$, and
$k(F)+r(F)=v(F)=v(G)$.  Therefore, $k(G)-k(F)+r(G)-r(F)=0$,
and we can rewrite our term as
$$A^{e(F)} B^{e(G)-e(F)} d^{\bc(F)-1}\ .$$
In terms of the state~$S = \vp^{-1}(F) \in \cS(\wt L_G)$
this summation term is equal to
$$A^{\a(S)} \, B^{\b(S)}\, d^{\de(S)-1}\,,$$
which is precisely the term of $\kb{\wt L_G}$ corresponding to the
state~$S \in \cS(\wt L_G)$.  This completes the proof.
\hspace{\fill}$\square$

\bigskip

\section{Final remarks and open problems.}\label{rem_kb}
\mbox{\ }\vspace{-15pt}

{\bf 1.} \,
Trivalent ribbon graphs are the main objects in
the finite type invariant theory of knots, links and
3-manifolds, while general ribbon graphs appeared in the
literature under a variety of different names
(see e.g.~\cite{DKC,BR2,Ka1}).
Embeddings of ribbon graphs into the 3-space are studied in~\cite{RT}.

\smallskip

{\bf 2.} \,
The Bollob\'as-Riordan polynomial can be defined by recurrent
contraction-deletion relations or by spanning tree expansion
similar to those of the Tutte polynomial, except that deletion by a loop
is not allowed.  We refer to~\cite{BR2,BR3} for the details.
We should note that~\cite{BR3} gives an extension to
unorientable surfaces as well. One can also find the contraction-deletion
relations and the spanning tree expansion for the signed
Bollob\'as-Riordan polynomial defined by (\ref{def_sbr}).
For a planar signed ribbon graph $\wh G$ the signed Bollob\'as-Riordan
polynomial $R_{\wh G}$ is related to
Kauffman's \emph{signed Tutte polynomial}~$Q[\wh G]$
from \cite{Ka2} by the formula
$$R_{\wh G}(x,y,z) \, = \, x^{\frac{v(\wh G)+1}{2}-k(\wh G)} \,
   y^{\frac{-v(\wh G)+1}{2}} \,
   Q[\wh G]\Bigl((y/x)^{1/2}, 1, (xy)^{1/2}\Bigr).
$$
So our version of $R_{\wh G}$ may be considered as a generalization
of the polynomial~$Q[\wh G]$ to signed ribbon graphs.
 If, besides the planarity, all edges of $\wh G$ are positive, then
 $R_{\wh G}$ is related to the dichromatic polynomial $Z[\Ga](q,v)$
 (see \cite{Ka2}) of the underlying graph $\Ga$:
 $$R_{\wh G}(x,y,z) \, = \, x^{-k(\wh G)} \,  y^{-v(\wh G)} \,
 Z[\Ga](xy,y)\, .
 $$

\smallskip

{\bf 3.} \,
It would be interesting to generalize the Bollob\'as-Riordan polynomial
for colored ribbon graphs \cite{BR1,Tr} and prove the corresponding
relation with the Kauffman bracket.  Let us also mention that
in~\cite{Ja} (see also~\cite{Tr}) Jaeger found a different relation
between links and graphs and proved that the whole Tutte polynomial,
not just its specialization, can be obtained from the HOMFLY
polynomial of the appropriate link.  Extending these results to
ribbon graphs is an important open problem.

Finally, recent results combinatorial evaluations of
the Tutte and Bollob\'as-Riordan polynomials~\cite{KP} leave an
open problem of finding such evaluations for general values of
polynomials~$R_G$.  It would be interesting to use the
Main Theorem to extend the results of~\cite{KP}.
We leave this direction to the reader.
\vskip.8cm

{\bf Acknowledgements}

\smallskip

\noindent
The second author is grateful to Jo Ellis-Monaghan,
Mike Korn and Vic Reiner for their insights into the
Tutte polynomial, and to B{\'e}la Bollob{\'a}s,
Vaughan Jones and Richard Stanley for encouragement.
The second author was supported by the NSA and the NSF.


\vskip1.cm

\parbox[t]{2.5in}{\it \textbf{Sergei~Chmutov}\\
Department of Mathematics\\
The Ohio State University, Mansfield\\
1680 University Drive\\
Mansfield, OH 44906\\
~\texttt{chmutov@math.ohio-state.edu}} \qquad\qquad
\parbox[t]{2.5in}{\it \textbf{Igor~Pak}\\
Department of Mathematics\\
Massachusetts Institute of Technology\\
77 Massachusetts Ave\\
Cambridge, MA 02139\\
~\texttt{pak@math.mit.edu}}

\end{document}